\documentclass[a4paper,12pt,makeidx]{amsart}
\usepackage{amsmath,amscd}
\usepackage{amssymb}
\usepackage{enumerate}
\usepackage{graphicx}
\usepackage{bbm}

\usepackage[noprefix]{nomencl}

\usepackage{hyperref} \hypersetup{colorlinks=true, bookmarks=false,
  urlcolor=black, linkcolor=black,citecolor=black}

\newtheorem*{Theorem}{Theorem}
\newtheorem*{Lemma}{Lemma}
\newtheorem*{Proposition}{Proposition}
\newtheorem*{Corollary}{Corollary}
\theoremstyle{definition} 
\newtheorem*{Remark}{Remark}

\newcommand{\Z}{\mathbbm{Z}}

\newcommand{\F}{\mathbbm{F}}
\renewcommand{\a}{\mathfrak{a}}

\newcommand{\gl}{\mathfrak{gl}}

\newcommand{\V}{V_\chi(\lambda)}
\renewcommand{\epsilon}{\varepsilon}
\renewcommand{\psi}{\varpsi}
\renewcommand{\:}{\colon}

\newcommand{\U}{U_{\chi}}

\DeclareMathOperator{\id}{id}

\DeclareMathOperator{\Soc}{Soc}

\DeclareMathOperator{\Hom}{Hom}
\DeclareMathOperator{\Ext}{Ext}
\DeclareMathOperator{\Img}{img}
\DeclareMathOperator{\End}{End}

\begin{document}
\title{Extensions of simple modules for the Witt algebra}
\author{Khalid Rian}\thanks{The author wishes to express his gratitude
  to Professor Jens Carsten Jantzen for his guidance and help.}
\address{Department of Mathematics, Aarhus University, Building 530,
  Ny Munkegade, 8000 Aarhus C, Denmark}
\email{khalid@imf.au.dk}
\begin{abstract}
The irreducible representations of the Witt algebra $W$ are completely
known. A classification of the irreducible $U_\chi(W)$--modules was first 
established by Chang and later simplified by Strade. The aim of this
article is to give a classification of the extensions of the simple
$U_\chi(W)$--modules having $p$--character $\chi$ of height $-1,0,1$ and $p-1$ where $p$ denotes the characteristic of the
ground field.
\end{abstract}
\maketitle
After completion of this paper, it has been brought to my attention
that similar results have appeared in \cite{BNW}. However, the methods
appearing in \cite{BNW} are different from those used here.
\section{Introduction}
\subsection{Notation}
Let $K$ be an algebraically closed field of characteristic $p>3$ and
let $W$ be the $p$--dimensional $\Z$--graded restricted simple Witt
algebra with basis $e_{-1},e_0,\dots,e_{p-2}$. The
Lie bracket is given by
\[
[e_i,e_j]=(j-i)e_{i+j} \quad \text{for all} \quad -1\leq i,j\leq p-2,
\]
where by definition $e_i=0$ for $i\not\in\{-1,0,\dots,p-2\}$. The
$p$--mapping is given by
\[
e_i^{[p]}=\delta_{i0}e_i\quad \text{for all} \quad -1\leq i\leq p-2.
\]
For each $ -1\leq i\leq p-1$ we set \nomenclature[110]{$W_i$}{\nomrefpage}
\[
W_i=Ke_i\oplus\cdots\oplus Ke_{p-2}
\]
and we define \textit{the height} $r(\chi)$ of a functional $\chi\in W^*$  to be
\[
r(\chi)=\operatorname{min}\{i\mid -1\leq i\leq p-1 \text{ and
}\chi_{|W_i}=0\}.
\]
\subsection{Simple modules} The simple $W$--modules were determined by
Chang in \cite{Chang}. We describe here very briefly the simple $W$--modules
with a $p$--character $\chi$ of height $-1$, $0$, $1$ and $p-1$.
\subsubsection{} Suppose that $\chi\in W^*$ has height at most 1. For
each weight $\lambda\in\Lambda(\chi)=\{\lambda\in K\mid
\lambda^p-\lambda=\chi(e_0)^p\}$ consider the $\chi-$reduced Verma module
$V_\chi(\lambda)=U_\chi(W)\otimes_{U_\chi(W_0)} K_\lambda$ where
$K_\lambda$ denotes $K$ considered as a $U_\chi(W_0)$--module via
\[
e_0\cdot1=\lambda\cdot1\quad\text{and}\quad e_i\cdot1=0 \ \text{  for all
 }\ i\neq0.
\]
\begin{Theorem}The following hold
\begin{enumerate}[(1)]
\item If $\chi$ has height $-1$, then there are $p$ isomorphism classes of simple $\U(W)$--modules. These modules are represented by the one-dimensional trivial $W$--module $K$, the $(p-1)$-dimensional module $S=V_\chi(p-1)/\Soc_{U_\chi(W)}V_\chi(p-1)$ and the $p$-dimensional modules $V_\chi(\lambda)$ for $\lambda\in\{1,2,\dots,p-2\}$. 
\item If $\chi$ has height $0$, then there are $p-1$ isomorphism classes of
simple $U_\chi(W)$--modules each of dimension $p$ and represented by
$V_\chi(\lambda)$ for $\lambda\in\{0,1,\dots,p-2\}$.
\item If $\chi$ has height $1$, then there are $p$ isomorphism classes of simple $U_\chi(W)$--modules each of dimension $p$ and represented by $V_\chi(\lambda)$ for $\lambda\in\Lambda(\chi)$.
\end{enumerate}
\end{Theorem}
For each $\mu\in\F_p=
\Z/p\Z$ we let $[\mu]\in\{0,1,\dots,p-1\}$ denote the unique
representative of $\mu$. Note that $\lambda-\lambda'\in\F_p$ for
all $\lambda,\lambda'\in\Lambda(\chi)$. Furthermore, if $\chi$ has height $1$ then $\lambda+\lambda'\not\in\F_p$ since otherwise this would imply $\lambda,\lambda'\in\F_p$.
\subsubsection{}If $\chi\in W^*$ has height $p-1$ then every simple
$U_\chi(W)$--module with one exception is projective. The remaining
simple module $L$ has a projective cover with two composition factors
both isomorphic to $L$, see \cite[Theorem 2.6]{FN}.
\begin{Remark}
If $1<r(\chi)<p-1$ then there exists a unique simple
$U_\chi(W)$--module up to isomorphism \cite[Hauptsatz 1]{Chang}. This module has a non-trivial self-extension. 
\end{Remark}
\section{Extensions of the reduced Verma modules}\label{Verma}
\subsection{}Let $\chi\in W^*$ with $r(\chi)=-1,0,1$ and let $\lambda,\lambda'\in\Lambda(\chi)$. The isomorphism
\[
\Ext_{U_\chi(W)}(V_\chi(\lambda'),V_\chi(\lambda))\simeq \Ext_{U_{\chi} (W_0)}(K_{\lambda'},V_{\chi} (\lambda))
\]
reduces the problem of classifying the extensions of Verma modules to
that of  classifying the extensions of $K_{\lambda'}$ by
$V_\chi(\lambda)$.
\subsection{}\label{Verma2}Consider a short exact sequence of $U_{\chi}(W_0)$--modules
\begin{equation}\label{seq1}
  0\rightarrow V_\chi(\lambda)\xrightarrow{f}M\xrightarrow{g} K_{\lambda'}\rightarrow 0.
\end{equation}
The Verma module $V_\chi(\lambda)$ has a basis $v_0,v_1,\dots,v_{p-1}$
given by $v_i=e_{-1}^i\otimes1$. Set $w_i=f(v_i)$ for all $i$. Since $U_\chi(Ke_0)$ is a semisimple algebra, there exists $w'\in M$ such that
$e_0w'=\lambda' w'$ and $g(w')=1$. Thus we obtain a basis for $M$
\begin{equation}\label{basis}
  w_0,w_1,\dots,w_{p-1},w',
\end{equation}
such that for all $0\leq j\leq p-2$
\begin{equation}\label{e_j1}
  e_jw_i=
  \begin{cases}
    (-1)^{j}\frac{i!}{(i-j)!}((j+1)\lambda-i+j)w_{i-j}, & \text{if }j\leq i,\\
    0, & \text{otherwise}.
  \end{cases}
\end{equation}
Furthermore, if $1\leq j\leq p-2$ then 
\[
e_0e_jw'=(e_je_0+je_j)w'=(\lambda'+j)e_jw',
\]
hence $e_jw'$ belongs to the 1-dimensional weight space $M_{\lambda'+j}$ spanned by $w_{[\lambda-\lambda'-j]}$. It follows that 
\begin{equation}
  \label{eq:3}
  e_jw'=\mathfrak{a}_jw_{[\lambda-\lambda'-j]} \quad \text{for some}\quad \mathfrak{a}_j\in K.
\end{equation}
\subsection{}\label{welldef}The $w'$ from the previous section is not unique. Since the $\lambda'$ weight space in $M$ is
spanned by $w'$ and $w_{[\lambda-\lambda']}$ any different
choice for $w'$ has the form $w'+bw_{[\lambda-\lambda']}$ for some
$b\in K$. Obviously, this leads to the same $e_jw'$ if
$e_jw_{[\lambda-\lambda']}=0$. Thus, if $\lambda=\lambda'$ then
$(\a_1,\a_2,\dots,\a_{p-2})$ is determined by the extension \eqref{seq1}. Furthermore, we have
\begin{Lemma}[A]
Let $\lambda,\lambda'\in\Lambda(\chi)$. If
$(\lambda,\lambda')=(0,p-1)$ or $(\lambda',\lambda)=(p-1,0)$ then $e_jw_{[\lambda-\lambda']}=0$ for all $j\geq1$.
\end{Lemma}
Note that the assumptions of the lemma imply $r(\chi)<1$ as $\lambda,\lambda'\not\in\F_p$ for $r(\chi)=1$. 
We let $\Theta(\chi)\subset \Lambda(\chi)\times\Lambda(\chi)$ denote the subset given by \nomenclature[230]{$\Theta(\chi)$}{\nomrefpage}
\[
(\{(0,p-1),(p-1,0)\}\cap\Lambda(\chi)\times\Lambda(\chi))\cup \{(\mu,\mu)\mid \mu\in\Lambda(\chi)\}.
\]
Lemma A and the discussion at the beginning of this section show that all the $\mathfrak{a_j}$ are determined by $M$ if $(\lambda,\lambda')\in\Theta(\chi)$. This is not the case for $(\lambda,\lambda')\not\in\Theta(\chi)$. However, we have
\begin{Lemma}[B] 
 Let $\lambda,\lambda'\in\Lambda(\chi)$ such that $(\lambda,\lambda')\not\in\Theta(\chi)$. If furthermore $\lambda+\lambda'\neq p-1$
  \textnormal{(}resp. $\lambda+\lambda'= p-1$\textnormal{)} then
  there is a unique choice for $w'$ such that $e_1w'=0$
  \textnormal{(}resp. $e_2w'=0$\textnormal{)}.  
\end{Lemma}
The proof is straightforward, given the fact that $e_1w_{[\lambda-\lambda']}\neq0$
(resp. $e_2w_{[\lambda-\lambda']}\neq0$) if $\lambda+\lambda'\neq p-1$ (resp. $\lambda+\lambda'= p-1$).
\begin{Remark}For each pair $(\lambda,\lambda')\not\in\Theta(\chi)$,
  we will always assume that $w'$ is the unique choice from the
  lemma. It follows, then, that the tuple $(\a_1,\a_2,\dots,\a_{p-2})$ determines
  the extension \eqref{seq1}. We call this tuple \textit{the $\a$--datum of \eqref{seq1}}. 
\end{Remark}
\subsection{}\label{Phi}Consider now a second extension
\begin{equation}\label{seq2}
  0\rightarrow V_\chi(\lambda)\xrightarrow{f'}M'\xrightarrow{g'}K_{\lambda'}\rightarrow 0
\end{equation}
of $\U(W_0)$--modules. If \eqref{seq1} and \eqref{seq2} are equivalent
extensions then there is an isomorphism $h:M\xrightarrow{\sim}M'$
compatible with the identities in both $V_\chi(\lambda)$ and $K_{\lambda'}$
\begin{equation*}
  \begin{CD}
    0 @>>> V_\chi(\lambda) @>f>> M @>g>> K_{\lambda'} @>>> 0\\
      @.      @|   @VVhV @| \\
    0 @>>> V_\chi(\lambda) @>f'>> M' @>g'>> K_{\lambda'} @>>> 0.
  \end{CD}
\end{equation*} 
We can choose $w''=h(w')$ to be the analogue to $w'$. Then $e_iw''=h(e_iw')=h(f(\mathfrak{a}_iv_{[\lambda-\lambda'-i]}))=f'(\mathfrak{a}_iv_{[\lambda-\lambda'-i]})$ for all $1\leq i\leq p-2$. Hence $f^{-1}(e_iw')$ depends only on the class of the extension \eqref{seq1} and we obtain a well-defined map
\begin{equation}\label{hom}
\Phi^i_{\lambda,\lambda'}:\Ext_{\U(W_0)}(K_{\lambda'},V_\chi(\lambda))\rightarrow K,
\end{equation}
that sends the class of \eqref{seq1} to $\a_i$.  One checks that
$\Phi^i_{\lambda,\lambda'}$ is linear. In Section \ref{a_i} we show
that every $\mathfrak{a}_i$ with $i\geq 3$ can be expressed as a
linear combination of $\mathfrak{a}_1$ and $\mathfrak{a}_2$.  It
follows that $\Phi_{\lambda,\lambda'} =(\Phi^1_{\lambda,\lambda'},\Phi^2_{\lambda,\lambda'})$ maps  $\Ext_{U_\chi(W_0)}(K_{\lambda'},V_\chi(\lambda))$ injectively into $K^2$. 
\begin{Remark}If $(\lambda,\lambda')\not\in\Theta(\chi)$ then  $\dim \Ext_{U_\chi(W_0)}(K_{\lambda'},V_\chi(\lambda))\leq1$. This is a consequence of Lemma \ref{welldef} B.
\end{Remark}
\subsection{}\label{a_i}Since $M$ is a $U_{\chi}(W_0)$--module, we have 
  \begin{equation*}
    [e_i,e_j]w'=(e_ie_j-e_je_i)w'  \ \text{ for all }i \text{ and }j.
  \end{equation*}
We will need the full strength of this formula later, but for now we
are content to remark that in the case where $i=1$ and $2\leq j\leq
p-3$, it yields
  \begin{align*}
    (j-1)\mathfrak{a}_{j+1}={}&-(-1)^j\biggl(\prod_{k=1}^{j}(\lambda-\lambda'-k)\biggr)(j(\lambda+1)+\lambda'+1)\mathfrak{a}_1\\
  &-(\lambda-\lambda'-j)(\lambda+\lambda'+j+1)\mathfrak{a}_j,
\end{align*}
or, equivalently, by induction
\begin{equation}\label{A1}
    \mathfrak{a}_j=A_j\mathfrak{a}_1+B_j\mathfrak{a}_2 \quad \text{for all}\quad 3\leq j\leq p-2,
  \end{equation}
  where
  \begin{align*}
    A_j={}&\frac{(-1)^j}{j-2}\biggl(\prod_{k=1}^{j-1}(\lambda-\lambda'-k)\biggr)\biggl((j-1)\lambda+\lambda'+j+\sum_{k=4}^{j}\frac{(j-k)!}{(j-3)!}\biggr.\\
    &\cdot((j+2-k)\lambda+\lambda'+(j+3-k))\prod_{l=0}^{k-4}(\lambda+\lambda'+j-l)\bigr),
  \end{align*}
  (the summation $\sum_{k=4}^j$ is understood to be 0 when $j=3$) and
  \[
  B_j=\frac{(-1)^j}{(j-2)!}\prod_{k=2}^{j-1}(\lambda-\lambda'-k)(\lambda+\lambda'+k+1).
  \]
Consequently, all the $\mathfrak{a}_i$ with $i\geq 3$ are determined by $\mathfrak{a}_1$ and $\mathfrak{a}_2$.
\subsection{}\label{M_a}Let $\lambda,\lambda'\in\Lambda(\chi)$. We
want to describe all possible $\a$--data of extensions as in
\eqref{seq1}. To this end, consider an arbitrary pair $(\a_1,\a_2)\in
K^2$ and use \eqref{A1} to extend it to a tuple
$\a=(\a_1,\a_2,\dots,\a_{p-2})$ in $K^{p-2}$. In order to simplify
notation set $\a_i=0$ for all $i>p-2$. Consider a vector space
$M_\a$ with a basis $w_0,w_1,\dots,w_{p-1},w'$ and define 
 endomorphisms $E_j\in\End_K(M_\mathfrak{a})$ for $0\leq j\leq p-2$ such that
\begin{equation*}
  E_jw_i=
  \begin{cases}
    (-1)^{j}\frac{i!}{(i-j)!}((j+1)\lambda-i+j)w_{i-j}, & \text{if }j\leq i,\\
    0, & \text{otherwise},
  \end{cases}
\end{equation*}
and
\begin{equation*}
  E_jw'=
  \begin{cases}
    \lambda'w', & \text{if }j=0,\\
    \mathfrak{a}_jw_{[\lambda-\lambda'-j]}, & \text{otherwise}.
  \end{cases}
\end{equation*}
Furthermore we set $E_j=0$ for $j\not\in\{0,1,\dots,p-2\}$. A comparison with the formulas in Section \ref{Verma2} shows that $\a$
is the $\a$-datum of an extension as in \eqref{seq1} if and only if
$M_\a$ has a $W_0$--module structure such that $e_jw=E_jw$ for all
$w\in M_\a$ and all $j\geq0$ and if this $W_0$--module then has
$p$--character $\chi$. We have on $M_\a$ a structure as $W_0$--module
if and only if 
\begin{equation}\label{act-1}
[E_i,E_j]w=(j-i)E_{i+j}w\quad\text{for all } i,j\geq0\text{ and } w\in M_\a.
\end{equation}
It suffices to check this for all $w$ in our basis. The linear map
$f\:\V\rightarrow M_\a$ with $f(v_i)=w_i$ for all $i$ satisfies
$E_jf(v)=f(E_jv)$ for all $v\in\V$ and all $j$. Consequently, \eqref{act-1}
holds for all $w_i$ and $M_\a$ has the desired structure if and only
if 
\begin{equation}\label{act0}
  [E_i,E_j]w'=(j-i)E_{i+j}w' \ \text{ for all }i,j\geq0.
\end{equation}
\subsection{}\label{W--module}
One verifies easily that \eqref{act0} holds for $i=0$. It also holds for
$i=j$. For $i=1$ and $2\leq j\leq p-3$ it is equivalent to \eqref{A1} which holds by the definition of the $\a_j$ ($j\geq3$). Furthermore, for every $(i,j)$ in 
\[
\{(1,p-2),(p-2,1)\}\cup\{(a,b)\in\Z^2\mid a\neq b \text{ and }2\leq a,b\leq p-2\},
\]
formula \eqref{act0} is equivalent to the following conditions 

\medskip

\noindent If $i>[\lambda-\lambda'-j]$ and $j>[\lambda-\lambda'-i]$, then
\begin{equation}\label{A2}
  \a_{i+j}=0.
\end{equation}
If $i>[\lambda-\lambda'-j]$ and $j\leq[\lambda-\lambda'-i]$, then 
\begin{equation}\label{A3}  (j-i)\a_{i+j}=(-1)^{j+1}\frac{[\lambda-\lambda'-i]!}{([\lambda-\lambda'-i]-j)!}(j(\lambda+1)+\lambda'+i)\a_i.
\end{equation}
If $i\leq[\lambda-\lambda'-j]$ and $j\leq[\lambda-\lambda'-i]$, then 
\begin{align}\label{A4}
  (j-i)\a_{i+j}={}& (-1)^{i}\frac{[\lambda-\lambda'-j]!}{([\lambda-\lambda'-j]-i)!}(i(\lambda+1)+\lambda'+j)\a_j\\
  &-(-1)^{j}\frac{[\lambda-\lambda'-i]!}{([\lambda-\lambda'-i]-j)!}(j(\lambda+1)+\lambda'+i)\a_i\notag.
\end{align}
\subsection{}\label{char}Suppose now that $\a$ satisfies \eqref{A2}--\eqref{A4} and hence
that $M_\a$ is a $W_0$--module with each $e_i$ acting as $E_i$. We
have then a short exact sequence 
\begin{equation}\label{seq11}
0\rightarrow\V\xrightarrow{f}M_\a\xrightarrow{g}K_{\lambda'}\rightarrow0.
\end{equation}
Evidently, $M_\a$ has $p$-character $\chi$ if and only if $e_i^p
w'=0$ for all $i\geq1$, or equivalently, if
$\a_ie_i^{p-1}w_{[\lambda-\lambda'-i]}=0$ for all $i\geq1$. We have
\[
\a_1e_1^{p-1}w_{p-1}=-\a_1\biggl(\prod_{j=0}^{p-2}(2\lambda-j)\biggr)w_0.
\] 
Observe that the term inside the bracket does not vanish if the height of $\chi$ is
1 because $\lambda\not\in\F_p$ in this case. Since furthermore
$e_i^{p-1}w_{[\lambda-\lambda'-i]}=0$ if $i>1$ or
$\lambda\neq\lambda'$ we obtain the following lemma
\begin{Lemma}$M_\a$ has $p$-character $\chi$ if and only if one of the
  following conditions holds
\begin{enumerate}
\item If $r(\chi)\in\{-1,0\}$ and $\lambda=\lambda'$ and $2\lambda=p-1$ then $\a_1=0$.
\item If $r(\chi)=1$ and $\lambda=\lambda'$ then $\a_1=0$.
\end{enumerate}
\end{Lemma}
\subsection{}We first look at the case where $\lambda=\lambda'$.
\begin{Proposition}\label{Self-Ext1}We have
  \begin{equation*}
    \Ext_{U_{\chi}(W_0)}(K_\lambda,V_{\chi}(\lambda))\simeq
    \begin{cases}
      K, &\text{if }\lambda\in\{0,p-1\} ,\\
      0, &\text{otherwise}.
    \end{cases}
  \end{equation*}
\end{Proposition}
It is convenient to break up the proof of the proposition into two lemmas. 
\begin{Lemma}[A]If $r(\chi)=1$ then
\[
\Ext_{U_\chi(W_0)}(K_\lambda,V_{\chi}(\lambda))=0.
\]
\end{Lemma}
\begin{proof}
Every extension of $K_\lambda$ by $V_{\chi}(\lambda)$ can be represented by a short exact sequence of $U_\chi(W_0)$--modules
\[
0\rightarrow V_{\chi}(\lambda)\rightarrow M_\a\rightarrow K_\lambda\rightarrow0.
\]
The claim to be proved amounts to saying that the above sequence splits, or equivalently, that $\a=0$. Since all the $\a_j$ can be expressed as a linear combination of $\a_1$ and $\a_2$, it suffices to prove that both $\a_1$ and $\a_2$ are 0. Note that $\a_1=0$ follows from Lemma \ref{char}. Thus, by definition, we have
\[
\a_j=(j-1)\biggl(\prod_{k=3}^{j}(2\lambda+k)\biggr)\a_2 \quad \text{for}\quad 3\leq j\leq p-2.
\]
Since the height of $\chi$ is 1 it follows that $2\lambda+k\neq0$ for all $k$. Therefore to prove $\a_2=0$ it suffices to show $\a_j=0$ for any $2\leq j\leq p-2$. Now, if we insert $(i,j)=(1,p-2)$ into \eqref{A4} and use the fact that $(p-1)!=-1$ we get
  \begin{equation*}
    -2(2\lambda-1)\a_{p-2}=0,
  \end{equation*}
proving the claim.  
\end{proof}
\begin{Lemma}[B]If $r(\chi)\in\{-1,0\}$ then
\[
\Ext_{U_\chi(W_0)}(K_\lambda,\V)\simeq
    \begin{cases}
      K, &\text{if }\lambda\in\{0,p-1\} ,\\
      0, &\text{otherwise}.
    \end{cases}
\]
\end{Lemma}
\begin{proof}
  Let $(\a_1,\a_2)$ be an arbitrary pair in $K^2$ and consider the vector space $M_\a$ constructed in Section \ref{M_a}. By definition, we have
  \begin{equation}\label{C1}
    \a_j=A_j\a_1+B_j\a_2 \quad \text{for all} \quad 3\leq j\leq p-2,
  \end{equation}
  where
  \begin{align*}    A_j=-\frac{\lambda+1}{j-2}\biggl(j!+\sum_{k=4}^{j}(k-1)!\prod_{l=k}^{j}\frac{l-1}{l-3}(2\lambda+l)\biggr),
  \end{align*}
  (the summation $\sum_{k=4}^j$ is understood to be 0 when $j=3$) and
  \begin{equation*}
    B_j=(j-1)\prod_{k=3}^{j}(2\lambda+k).
  \end{equation*}
  We rewrite \eqref{A2}--\eqref{A4} for the present case. First, note that \eqref{A2} and \eqref{A3} hold trivially since $[\lambda-\lambda-i]=p-i$ and $[\lambda-\lambda-j]=p-j$ for all $1\leq i,j\leq p-2$. If we insert $(i,j)=(1,p-2)$ into \eqref{A4} and use $(p-1)!=-1$, we
  see that
  \begin{equation}\label{C2}
    2(2\lambda-1)\a_{p-2}=(\lambda+1)\a_1.
  \end{equation}
 If $2\leq i,j\leq p-2$ and  $i+j\leq p$, the same formula yields 
  \begin{align}\label{C3}
    (j-i)\a_{i+j}={}&\frac{(i+j-1)!}{(j-1)!}((i+1)\lambda+i+j)\a_j\\  
    &-\frac{(i+j-1)!}{(i-1)!}((j+1)\lambda+i+j)\a_i\notag. 
  \end{align}

  We now proceed by showing how these formulas are related. Let's
  begin by assuming that $\lambda\not\in\{0,p-1\}$. If $2\lambda=p-1$ then Lemma \ref{char} implies $\a_1=0$, which together with \eqref{C1} implies
  \begin{equation}\label{A1a}
    \a_j=(j-1)(j-1)!\a_2 \quad \text{for all}\quad 3\leq j\leq p-2.
  \end{equation}
It follows that $\a_2=0$ because otherwise this would mean that $\a_{p-2}\neq0$ and then by \eqref{C2} that $2\lambda=1$  in contradiction with our assumption.

\medskip

Suppose next that $2\lambda\neq p-1$.  If we insert $(i,j)=(2,p-2)$ into \eqref{C3}, we obtain 
  \begin{equation}\label{p-2}
    6\a_{p-2}=\a_2.
  \end{equation}

\medskip

\noindent Assume first that $2\lambda=1$. Since $\lambda\neq p-1$ eq.  \eqref{C2} implies $\a_1=0$ hence \eqref{C1} reduces to 
  \begin{equation*}
    \a_j=\frac{1}{6}(j-1)(j+1)!\a_2 \quad \text{for all} \quad 3\leq j\leq p-2.
  \end{equation*} 
  In  particular, we have $2\a_{p-2}=\a_2$, which together with \eqref{p-2}, implies $\a_2=0$, as desired. 

\medskip

\noindent Assume next that $2\lambda\neq1$. Eq. \eqref{C2} and \eqref{p-2}
  give the following relation
  \begin{equation}\label{rel:a_1,a_2}
    3(\lambda+1)\a_1=(2\lambda-1)\a_2,
  \end{equation}
  which means that $\a_1$ and $\a_2$ are either both 0 or both nonzero. We assume the latter and derive a contradiction. By definition, we have
  $\a_3=-6(\lambda+1)\a_1+2(2\lambda+3)\a_2$ hence it follows from \eqref{rel:a_1,a_2} that $\a_3=8\a_2$.
  Thus, if we insert $(i,j)=(3,p-3)$ into \eqref{C3}, we obtain
  \begin{equation*}
    3\a_{p-3}=-\a_2.
  \end{equation*} 
Likewise, if we  insert $(i,j)=(2,p-3)$ into \eqref{C3}, we obtain
  \[
  6(3\lambda-1)\a_{p-3}=-(2\lambda+1)\a_2.
  \]
  Since $2\lambda+1\neq0$ we can determine $\lambda$ by eliminating $\a_2$ and $\a_{p-3}$
  \[
  4\lambda=3.
  \]
Note $3\a_{p-3}=-\a_2$ implies $p>5$. If we insert $(i,j)=(2,3)$ into \eqref{C3}, we see 
  \[
 \a_5=(192\lambda+360)\a_2=504\a_2.
  \]
However, by definition we have
\[
\a_5=\frac{1}{3}(128\lambda^2+656\lambda+1080)\a_2=548\a_2,
\]
which implies first that $p=11$ and then that $\lambda=9$.  If we put
all this together and use \eqref{C1} we see $\a_1=\a_4$. But inserting
$(i,j)=(4,6)$ into \eqref{C3} implies $\a_4=0$ and thus $\a_1=0$. Contradiction!

\medskip

We now move to the case where $\lambda=0$. A necessary condition for $M_\a$ to be a $U_\chi(W_0)$--module is that $2\a_{p-2}=-\a_1$, cf. \eqref{C2}. Since by definition 
\[
\a_{p-2}=4\a_1-\frac{3}{2}\a_2,
\]
this implies $\a_2=3\a_1$.  Now, the claim to be proved amounts to saying that each pair $(\a_1,\a_2)$ with $\a_2=3\a_1$ induces a $U_\chi(W_0)$--module $M_\a$ in the way described in Section \ref{M_a}. Using our previous notation, we have
  \[
  A_j=-(j-2)j! \quad  \text{and} \quad  B_j=\frac{1}{2}(j-1)j!,
  \]
  hence $\a_2=3\a_1$ implies
  \begin{equation*}
    \a_j=\frac{1}{2}(j+1)!\a_1 \quad \text{for all} \quad 1\leq j\leq p-2.
  \end{equation*}
  In particular, we have $\a_{p-2}=-1/2\a_1$ in consistence with \eqref{C2}. A straightforward computation shows that
  \begin{align*}
    (j-i)A_{i+j}&=(i+j)!\left(\frac{A_j}{(j-1)!}-\frac{A_i}{(i-1)!}\right),\\
    (j-i)B_{i+j}&=(i+j)!\left(\frac{B_j}{(j-1)!}-\frac{B_i}{(i-1)!}\right).
  \end{align*}
  proving \eqref{C3} for $i+j\leq p-2$. Likewise, a simple computation shows 
  \[
  \a_j(p-j-2)!=\a_{p-j-1}(j-1)!,
  \] 
proving \eqref{C3} for $i+j=p-1$. Since $i+j=p$ gives 0 on both sides of \eqref{C3} the claim follows from Lemma \ref{char}.

\medskip 

Finally, we consider the case where $\lambda=p-1$. Here we have $A_j=0$ and it
follows from \eqref{C2} that $\a_{p-2}=0$. Therefore since 
 \[
 \a_j=(j-1)!\a_2 \quad \text{for all}  \quad 2\leq j\leq p-2,
 \]
this implies first $\a_2=0$ and then that $\a_j=0$ for all $j\geq2$. The rest of the proof is straightforward.
\end{proof}
\subsection{}\label{(0,p-1)} Next, we address the cases where $(\lambda,\lambda')=(0,p-1)$ and $(\lambda,\lambda')=(p-1,0)$. This can only occur if the height of $\chi$ is less than 1 because $\lambda,\lambda'\not\in\F_p$ for $r(\chi)=1$. The arguments presented previously would apply equally well to these cases, but we will give another, simpler, proof which relies on the fact that $\Ext_{U_\chi(W_0)}(K_0,V_\chi(0))$ is 1-dimensional, cf. Proposition \ref{Self-Ext1}.
\begin{Proposition} If $(\lambda,\lambda')=(0,p-1)$ or $(\lambda,\lambda')=(p-1,0)$, then 
  \begin{equation*}
    \Ext_{U_\chi(W_0)}(K_{\lambda'},V_\chi(\lambda))\simeq  K.
  \end{equation*}
\end{Proposition}
\begin{proof}
Let $\eta\in W^*$ be a functional of height $0$. We have
\[
\Ext_{U_\eta(W)}(V_\eta(0),V_\eta(0))\simeq\Ext_{U_\chi(W_0)}(K_0,V_\chi(0))\simeq  K.
\] 
Since the reduced Verma modules $V_\eta(0)$ and $V_\eta(p-1)$ are isomorphic \cite [Hilfssatz 7]{Chang}, we obtain
\begin{equation*}
  \Ext_{U_\chi(W_0)}(K_0,V_\chi(p-1))\simeq\Ext_{U_{\eta}(W)}(V_\eta(0),V_\eta(p-1))\simeq K,
\end{equation*}
as claimed. The other case can be handled similarly. 
\end{proof}
\subsection{}\label{neq p-1}It is convenient to divide the case $(\lambda,\lambda')\not\in \Theta(\chi)$ into two subcases depending on whether or not $\lambda+\lambda'=p-1$. Recall that $\lambda+\lambda'=p-1$ can only occur if $r(\chi)<1$. The following lemma will be needed in the sequel.
\begin{Lemma}
Let $i$ and $j$ be two integers with $2\leq i,j\leq p-2$. 
\begin{enumerate}
\item If $i>[\lambda-\lambda'-j]$ then $i+j>[\lambda-\lambda']$. 
\item If $j\leq[\lambda-\lambda'-i]$ and $i+j>[\lambda-\lambda']$ then $i>[\lambda-\lambda']$.
\end{enumerate}
\end{Lemma}
\subsection{}We denote the roots of the polynomial $x^2-a\in K[x]$ by $\pm\sqrt{a}$. \nomenclature[250]{$\pm\sqrt{a}$}{\nomrefpage}
\begin{Proposition}\label{neq}If $(\lambda,\lambda')\not\in \Theta(\chi)$ and $\lambda+\lambda'\neq p-1$, then
\begin{equation*}
\begin{array}{c}
\Ext_{U_{\chi}(W_0)}(K_{\lambda'},V_{\chi}(\lambda))\\ \ \\
 \simeq\begin{cases}
      K,&\text{if }[\lambda-\lambda']=1 \text{ and }\lambda\in\{1,p-1\},\\
      K,&\text{if }[\lambda-\lambda']\in\{2,3\},\\
      K,&\text{if }[\lambda-\lambda']=4 \text{ and }p>5,\\   
      K,&\text{if }[\lambda-\lambda']=4 \text{ and }\lambda\in\{0,3\},\\
      K,&\text{if }[\lambda-\lambda']= 5 \text{ and }\lambda\in\{0,4\},\\ 
      K,&\text{if }[\lambda-\lambda']= 6 \text{ and
      }2\lambda=5\pm\sqrt{19}\text{ and }p>7,\\
      0,&\text{otherwise}.
    \end{cases}
\end{array}
\end{equation*}
\end{Proposition}
\begin{proof}
Every extension of $K_{\lambda'}$ by $V_{\chi}(\lambda)$ can be represented by a short exact sequence of $U_{\chi}(W_0)$--modules
  \begin{equation}\label{seqv}
    0\rightarrow V_{\chi}(\lambda)\rightarrow M_\a\rightarrow K_{\lambda'}\rightarrow 0,
  \end{equation} 
where $\a_1=0$ (see Lemma \ref{welldef} B and Remark \ref{welldef}) and
  \begin{equation}\label{D1}
    \a_j=\frac{(-1)^j}{(j-2)!}\biggl(\prod_{k=2}^{j-1}(\lambda-\lambda'-k)(\lambda+\lambda'+k+1)\biggr)\a_2, \quad 3\leq j\leq p-2.
  \end{equation}
Note that $\a_j=0$ for all $j>[\lambda-\lambda']$ if $[\lambda-\lambda']\geq2$. This is a very useful observation which we shall use several times, often without any reference. We divide the proof into several steps depending on $[\lambda-\lambda']$.

 \medskip 

\noindent\textbf{Case 1.} Suppose that  $[\lambda-\lambda']=1$.

 \medskip

\noindent  For the sake of simplicity, we assume that $p>5$. Eq. \eqref{D1} becomes
\[
\a_j=\biggl(\prod_{k=2}^{j-1}(2\lambda+k)\biggr)\a_2 \quad \text{for all} \quad 3\leq j\leq p-2.
\] 
If we insert $(i,j)=(2,3)$ into \eqref{A4} and then use the above
formula to express $\a_3$ as
a scalar multiple $\a_2$ we obtain
  \[
  \a_5=(2\lambda+2)(18\lambda+12)\a_2.
  \]
A necessary condition for \eqref{seqv} to be non-split is that $\a\neq0$, or equivalently, that $\a_2\neq0$. If this is the case, then
\[
\prod_{k=2}^{4}(2\lambda+k)=(2\lambda+2)(18\lambda+12),
\]
which implies $\lambda\in\{1,p-1\}$. (It should be noted that $\lambda=0$ is a solution to the above equation, but it has been deliberately omitted because otherwise we would have $\lambda'=p-1$ and thus $\lambda+\lambda'=p-1$.) 

\medskip

Conversely, suppose that $\lambda=1$ or $\lambda=p-1$. Since
$\lambda-\lambda'=1$, this implies $\lambda'=0$ or $\lambda'=p-2$,
respectively. Let $(\a_1,\a_2)$ be a pair in $K^2$ such that $\a_1=0$
and consider the corresponding vector space $M_\a$. If $\lambda=p-1$,
then $\lambda+\lambda'=p-3$ hence by definition  $\a_j=0$ for all
$j\neq2$ and a very simple computation shows that \eqref{A2}--\eqref{A4} hold. Suppose next that $\lambda=1$. We have then
\[
p-2>[\lambda-\lambda'-1]\quad\text{and}\quad 1<[\lambda-\lambda'-(p-2)].
\]
If we insert $(i,j)=(p-2,1)$ into \eqref{A3}, we obtain  zero on both sides because $(\lambda+1)+\lambda'-2=0$ and $\a_{p-1}=0$. Furthermore, if $2\leq i,j\leq p-2$ then
\[
[\lambda-\lambda'-i]=p+1-i \quad \text{and}\quad [\lambda-\lambda'-j]=p+1-j.
\]
Thus, we may assume that 
\[
i\leq[\lambda-\lambda'-j] \quad \text{and}\quad j\leq[\lambda-\lambda'-i],
\]
or equivalently that $i+j\leq p+1$. We have
\[
\a_k=\frac{1}{6}(k+1)!\a_2 \quad \text{for all}\quad 2\leq k\leq p-2,
\]
which, when inserted into \eqref{A4}, yields
\begin{align*}
(j-i)\a_{i+j}={}& \frac{1}{6}(i+j-2)!(2i+j)(j-1)j(j+1)\a_2\\
    &-\frac{1}{6}(i+j-2)!(2j+i)(i-1)i(i+1)\a_2.
\end{align*}
Since
\[
(j-i)\a_{i+j}=
\begin{cases}
\frac{1}{6}(j-i)(i+j+1)!\a_2, &\text{if }i+j\leq p-2,\\
0, &\text{if }i+j\in\{p,p\pm1\},
\end{cases}
\]
the claim follows from a straightforward computation.

\medskip

\noindent\textbf{Case 2.}  Suppose that $[\lambda-\lambda']\in\{2,3\}$. 

\medskip

\noindent Let $(\a_1,\a_2)$ be a pair in $K^2$ such that $\a_1=0$ and
consider the corresponding vector space $M_\a$. Lemma \ref{neq p-1}
together with the remark following eq. \eqref{D1} yields \eqref{A2}--\eqref{A4}
for $2\leq i,j\leq p-2$. Since furthermore $\a_{p-2}=0$ for $p>5$ this
proves \eqref{A2}--\eqref{A4} for all $i,j$ in the case $p>5$. (Note that the above reasoning works just as well for $\lambda-\lambda'=4$.) The case $p=5$ is straightforward and is omitted here.
\medskip

\noindent\textbf{Case 3.} Suppose that $[\lambda-\lambda']=4$.

\medskip

\noindent As noted above, we may assume $p=5$ and $(i,j)\in\{(1,3),(3,1)\}$. We have
\[
3=[\lambda-\lambda'-1] \quad \text{and}\quad 1=[\lambda-\lambda'-3].
\]
Thus, if we insert $(i,j)=(1,3)$ into \eqref{A4} and then use the definition of $\a_3$ to express it as a scalar multiple $\a_2$, we obtain
\[
2(\lambda+\lambda'+3)(\lambda+\lambda'+4)\a_2=0.
\]
Now, every extension of $K_{\lambda'}$ by $V_{\chi}(\lambda)$ can be represented by a short exact sequence as in \eqref{seqv} and such that $\a_1=0$. The above computation shows that a necessary (and sufficient) condition for \eqref{seqv} to be non-split is that $\lambda+\lambda'\in\{1,2\}$, or equivalently, that $\lambda\in\{0,3\}$.

\medskip

\noindent\textbf{Case 4.} Suppose that $[\lambda-\lambda']=5$.

\medskip

\noindent First, note that the assumption implies $p>5$ since $[\lambda-\lambda']<p$. We consider a short exact sequence of $U_{\chi}(W_0)$--modules as in \eqref{seqv}. If we insert
  $(i,j)=(2,3)$ into \eqref{A4} and then use the definition of
  $\a_5$ to express it as a scalar multiple $\a_2$, we get 
  \begin{align}\label{23} \frac{1}{6}\biggl(\prod_{k=3}^5(\lambda+\lambda'+k)\biggr)\a_2=&\ (2\lambda+\lambda'+5)(\lambda+\lambda'+3)\a_2\\ \notag
&-(3\lambda+\lambda'+5)\a_2,
  \end{align}
  which since $\lambda+\lambda'\neq p-1$ can be rewritten as
  \begin{equation*}
((\lambda+\lambda')(\lambda+\lambda'+5)-6\lambda)\a_2=0.
  \end{equation*}
Thus, the term inside the parentheses must equal 0 in order for \eqref{seqv} to be non-split. Note that the above reasoning works just as well in the case $[\lambda-\lambda']\geq5$. We shall make use of this later, but for now we are 
content to remark that in our present case this implies $\lambda\in\{0,4\}$.

\medskip

Conversely, suppose that $\lambda\in\{0,4\}$ and let $(\a_1,\a_2)$ be
a pair in $K^2$ such that $\a_1=0$. We extend $(\a_1,\a_2)$ as usual
to a tuple $\a\in K^{p-2}$. The choice of $\lambda$ together with
Lemma \ref{neq p-1} implies \eqref{A2}--\eqref{A4} for $2\leq i,j\leq
p-2$ hence we may assume $(i,j)\in\{(1,p-2),(p-2,1)\}$.

If $p=7$ then
\[
p-2>[\lambda-\lambda'-1]\quad \text{and}\quad 1>[\lambda-\lambda'-(p-2)].
\]
Since by definition $\a_{p-1}=0$, this proves \eqref{A2}. If $p>7$, then
\[
p-2>[\lambda-\lambda'-1]\quad \text{and}\quad 1<[\lambda-\lambda'-(p-2)].
\]
The remark following eq. \eqref{D1} implies $\a_{p-2}=0$ and thus \eqref{A3}, thereby proving the claim.

\medskip

\noindent\textbf{Case 5.} Suppose that $[\lambda-\lambda']=6$.

\medskip

\noindent As previously noted, every nontrivial extension of
$K_{\lambda'}$ by $V_{\chi}(\lambda)$ can be represented by a sequence
of $U_{\chi}(W_0)$--modules as in \eqref{seqv} and such that $\a_1=0$
and $(\lambda+\lambda')(\lambda+\lambda'+5)=6\lambda$. In our present case, this amounts to
\[
  \lambda=\frac{1}{2}(5\pm\sqrt{19}).
\]
For $r(\chi)<1$ we have $\lambda\in\F_p$ and since $x^2-19$ does not split in $\F_7[x]$ this shows that $p>7$. For $r(\chi)=1$ we have $\lambda+\lambda'\not\in\F_p$. If $p=7$ then
\[
p-2=[\lambda-\lambda'-1] \quad \text{and}\quad 1=[\lambda-\lambda'-(p-2)].
\]
But if we insert $(i,j)=(1,p-2)$ into \eqref{A4} and then express $\a_{p-2}$ as a multiple of $\a_2$ we obtain
\[
\frac{1}{6}(\lambda+\lambda'-1)\biggl(\prod_{k=2}^{4}(\lambda-\lambda'-k)(\lambda+\lambda'+k+1)\biggr)\a_2=0.
\]
This is only possible if $\a_2=0$, or equivalently, if $\a=0$ contradicting the fact that the sequence \eqref{seqv} does not split.

\medskip

Conversely, assume that $p>7$ and let $\lambda$ be as above. Let
$(\a_1,\a_2)\in K^2$ be a pair such that $\a_1=0$ and consider the
corresponding vector space $M_\a$. The remark following eq. \eqref{D1}
gives \eqref{A2}--\eqref{A4} for $(i,j)\in\{(1,p-2),(p-2,1)\}$. By
virtue of Lemma \ref{neq p-1} we may assume that
$i+j\leq[\lambda-\lambda']$. Since furthermore we have chosen $\lambda$ in such a way that \eqref{A4} holds for $(i,j)=(2,3)$, it suffices to consider the case $(i,j)\in\{(2,4),(4,2)\}$. We have
\[
4=[\lambda-\lambda'-2] \quad \text{and}\quad 2=[\lambda-\lambda'-4].
\]
If we insert $(i,j)=(2,4)$ into \eqref{A4} and then use the definition of $\a_4$ to express it as a scalar multiple $\a_2$, we obtain
  \begin{align*}
\frac{1}{12}\biggl(\prod_{k=3}^6(\lambda+\lambda'+k)\biggr)\a_2={}&  \frac{1}{2}(2\lambda+\lambda'+6)\biggl(\prod_{k=3}^4(\lambda+\lambda'+k)\biggr)\a_2\\ \notag
&-(4\lambda+\lambda'+6)\a_2,
  \end{align*}
which gives the equation
\[ 
(\lambda+\lambda'+1)(\lambda+\lambda'+6)((\lambda+\lambda')(\lambda+\lambda'+5)-6\lambda)\a_2=0,
\]
that holds by the definition of $\lambda$.

\medskip

\noindent\textbf{Case 6.} Suppose that $[\lambda-\lambda']\geq7$.

\medskip

\noindent Every extension of $K_{\lambda'}$ by $V_{\chi}(\lambda)$ can
be represented by a short exact sequence of $U_{\chi}(W_0)$--modules
as in \eqref{seqv} and such that $\a_1=0$. A necessary condition for
\eqref{seqv} to be non-split is that
$(\lambda+\lambda')(\lambda+\lambda'+5)=6\lambda$ which, in our present case, implies
\begin{equation}\label{relat}
  \lambda'=\frac{1}{2}( -5\pm\sqrt{24\lambda+25})-\lambda.
\end{equation}
Assume towards contradiction that $\a_2\neq0$. We have
\[
5\leq[\lambda-\lambda'-2] \quad \text{and}\quad 2\leq[\lambda-\lambda'-5].
\]
Thus, if we insert $(i,j)=(2,5)$ into \eqref{A4} and then use the definition of $\a_5$ to write it as a scalar multiple $\a_2$, we get
  \begin{align*}
 \frac{1}{40}\biggl(\prod_{k=3}^7(\lambda+\lambda'+k)\biggr)\a_2={}& \frac{1}{6}(2\lambda+\lambda'+7)\biggl(\prod_{k=3}^5(\lambda+\lambda'+k)\biggr)\a_2\\ \notag
&-(5\lambda+\lambda'+7)\a_2.
  \end{align*}
which, together with \eqref{relat}, yields
\[
\pm(\lambda^2+5\lambda+4)\sqrt{24\lambda+25}=15\lambda^2+35\lambda+20.
\]
Squaring both sides and subtracting one from the other yields the equation
\[
\frac{9}{50}\lambda(\lambda-1)(\lambda+1)^2(3\lambda+2)=0,
\]
which in turn yields $\lambda\in\{0,\pm1,-2/3\}$. Now we can insert this back into \eqref{relat} to determine $\lambda'$. But first note that if
\begin{equation*}
  \lambda'=\frac{1}{2}( -5-\sqrt{24\lambda+25})-\lambda,
\end{equation*} 
then $\lambda=p-1$ because $0,1,-2/3$ are not solutions to
\begin{equation*}
-(\lambda^2+5\lambda+4)\sqrt{24\lambda+25}=15\lambda^2+35\lambda+20.
\end{equation*}
Putting all of this together, we obtain
\[
(\lambda,\lambda')\in\{(p-1,p-1),(p-1,p-2),(1,0),(0,0),(-2/3,-1/3)\}.
\] 
But none  of these cases is possible since by assumption 
  $[\lambda-\lambda']\geq7$ and  $\lambda+\lambda'\not=p-1$; a contradiction which can only be avoided if $\a_2=0$. This completes the proof of the proposition.
\end{proof}
\subsection{}We are left with the case where
$(\lambda,\lambda')\not\in \Theta(\chi)$ and
$\lambda+\lambda'=p-1$. Note that this implies $r(\chi)<1$.
\begin{Proposition} If $(\lambda,\lambda')\not\in \Theta(\chi)$ and $\lambda+\lambda'=p-1$, then 
  \begin{equation*}
    \Ext_{U_\chi(W_0)}(K_{\lambda'},\V)\simeq  
\begin{cases}
 K,&\text{if } [\lambda-\lambda']\in\{2,3,4\}, \\
 K,&\text{if } [\lambda-\lambda']=6 \text{ and }p=19,\\
 0,&otherwise.
\end{cases}
  \end{equation*}
\begin{proof}
  Let $(\lambda,\lambda')\not\in \Theta(\chi)$ such that
  $\lambda+\lambda'= p-1$. Every extension of $K_{\lambda'}$ by
  $V_\chi(\lambda)$ can be represented by a short exact sequence of
  $U_\chi(W_0)$--modules as \eqref{seqv}. The assumption
  $\lambda+\lambda'= p-1$ implies $\a_2=0$, see Lemma \ref{welldef} B
  and Remark \ref{welldef}. Using the foregoing notation, we have
  \begin{align*} A_j={}&\frac{(-1)^j}{j-2}\biggl(\prod_{k=0}^{j-2}(2\lambda-k)\biggr)\biggl((j-2)\lambda-1+j+\sum_{k=4}^{j}\frac{(j-k)!}{(j-3)!}\biggr.\\ \notag
    &\cdot((j+1-k)\lambda+(j+2-k))\prod_{l=1}^{k-3}(j-l)\biggr).
  \end{align*}
(The summation $\sum_{k=4}^j$ is understood to be 0  when $j=3$.) Note that $\a_j=0$ for all $j>[\lambda-\lambda']$. Furthermore, it should also be noted that $\lambda-\lambda'\neq1$ since otherwise we would have $(\lambda,\lambda')=(0,p-1)$ in contradiction to the assumption that $(\lambda,\lambda')\not\in\Theta(\chi)$.

\medskip 

\noindent\textbf{Case 1.} Suppose that $[\lambda-\lambda']\in\{2,3,4\}$.

\medskip

\noindent Conversely, every pair $(\a_1,\a_2)\in K^2$ with $\a_2=0$
gives rise to a module $M_\a$ as described in Section \ref{M_a}: Lemma
\ref{neq p-1} yields \eqref{A2}--\eqref{A4} for $2\leq i,j\leq
p-2$. For $(i,j)\in\{(1,p-2),(p-2,1)\}$ the computations are simple
and are omitted for brevity.

\medskip 

\noindent\textbf{Case 2.} Suppose that $[\lambda-\lambda']=5$.

\medskip 

\noindent Consider a short exact sequence of $U_\chi(W_0)$--modules as in \eqref{seqv}. If we insert $(i,j)=(2,3)$ into \eqref{A4} and then insert the expressions $\a_3=A_3\a_1$ and $\a_5=A_5\a_1$ into the result we get
  \begin{align*}\lambda(\lambda+2)(\lambda+4)\biggl(\prod_{k=1}^3(2\lambda-k)\biggr)\a_1=\frac{1}{3}\lambda(13\lambda+22)\biggl(\prod_{k=1}^3(2\lambda-k)\biggr)\a_1.
  \end{align*}
 Now, the assumptions $\lambda+\lambda'=p-1$ and $(\lambda,\lambda')\not\in\Theta(\chi)$ imply $\lambda\neq0$. Since, in addition, $[\lambda-\lambda']=5$, the above equation reduces to
  \[
  3(\lambda+2)(\lambda+4)\a_1=(13\lambda+22)\a_1,
  \]
  which, by subtracting the right side from the left side, gives
  \[
  (3\lambda^2+5\lambda+2)\a_1=0.
  \]
We get a quadratic equation which can be solved easily; the roots are $p-1$ and $-2/3$. But, we  can immediately exclude $\lambda=p-1$ since this would imply $(\lambda,\lambda')=(p-1,0)$. All the above reasoning works just as well in the case where $[\lambda-\lambda']\geq5$. We shall make use of this later, but for now we are content to remark that in our present case $\lambda+\lambda'=p-1$ implies $(\lambda,\lambda')=(-2/3,-1/3)$ contradicting the fact that $\lambda-\lambda'=5$.

\medskip 

\noindent\textbf{Case 3.} Suppose that $[\lambda-\lambda']=6$.

\medskip 

\noindent As noted previously, a necessary condition for \eqref{seqv}
to be non-split is that $\lambda=-2/3$. The assumption
$\lambda+\lambda'=p-1$ implies $\lambda-\lambda'=2\lambda+1$ which
since $\lambda-\lambda'=6$ means $\lambda=5/2$ and hence
$p=19$. Conversely, let $p=19$ and $(\lambda,\lambda')=(12,6)$. (Note
$12=5/2$ in characteristic $p=19$.) Let $(\a_1,\a_2)$ be a pair in
$K^2$ such that $\a_2=0$ and consider the corresponding vector space
$M_\a$. Eq. \eqref{A2}--\eqref{A4} follow immediately for
$(i,j)\in\{(1,p-2),(p-2,1)\}$ since $\a_{p-2}=0$. By virtue of Lemma \ref{neq p-1} we may assume that $(i,j)\in\{(2,4),(4,2)\}$. We have
\[
2=[\lambda-\lambda'-4]\quad \text{and}\quad 4=[\lambda-\lambda'-2].
\]
If we insert $(i,j)=(2,4)$ into \eqref{A4} and then insert the expressions $\a_4=A_4\a_1$ and $\a_6=A_6\a_1$ into the result, we get
\[ \frac{1}{3}\lambda(77\lambda+125)\biggl(\prod_{k=1}^{4}(2\lambda-k)\biggr)\a_1=\lambda(\lambda+5)(5\lambda+9)\biggl(\prod_{k=1}^{4}(2\lambda-k)\biggr)\a_1.
\]
Keeping in mind that the characteristic is 19, a very simple
computation shows that the above equation holds for all $\a_1$.

\medskip 

\noindent\textbf{Case 4.} Suppose that $[\lambda-\lambda']\geq7$.

\medskip

\noindent Suppose we are given a short exact sequence of
$U_\chi(W_0)$--modules as in \eqref{seqv} and such that $\a_2=0$ and
$\lambda=-2/3$. We have
\[
2\leq[\lambda-\lambda'-5]\quad\text{and}\quad 5\leq[\lambda-\lambda'-2].
\]
If we insert $(i,j)=(2,5)$ into \eqref{A4}, we obtain
  \begin{align*} \frac{6}{10}\lambda(87\lambda+137)\biggl(\prod_{k=1}^{5}(2\lambda-k)\biggr)\a_1={}&\frac{2}{3}\lambda(\lambda+6)(13\lambda+22)\\
&\cdot\biggl(\prod_{k=1}^{5}(2\lambda-k)\biggr)\a_1.
  \end{align*}
Since $\lambda\neq0$ and $[\lambda-\lambda']\geq7$, the above equation reduces to
\[
9\left(87\lambda+137\right)\a_1=10(\lambda+6)(13\lambda+22)\a_1,
\]
 which, by subtracting the left side from the right side, gives 
\[
(130\lambda^2+217\lambda+87)\a_1=0.
\]
The claim follows since $-2/3$ is not a root of the polynomial inside the brackets. 
\end{proof}
\end{Proposition}
\subsection{}\label{Theorem1}We summarize the preceding results as a theorem as follows.
\begin{Theorem}\label{sum1} We have the following three cases
\begin{enumerate}
\item If $\lambda'\in\{0,p-1\}$, then
  \begin{equation*}
    \Ext_{U_\chi(W)}(V_\chi(\lambda'),V_\chi (\lambda))\simeq
    \begin{cases}
      K,& \text{if }\lambda\in\{0,1,2,3,4,p-1\},\\
      0,&\text{otherwise}.
    \end{cases}
  \end{equation*}
\item If $\lambda\in\{0,p-1\}$, then
  \begin{equation*}
    \Ext_{U_\chi(W)}(V_\chi(\lambda'),V_\chi(\lambda))\simeq
    \begin{cases}
      K,& \text{if }\lambda'\in\{0\}\cup\{p-i\mid 1\leq i \leq 5\},\\
      0,&\text{otherwise}.
    \end{cases}
  \end{equation*}
\item If $\lambda,\lambda'\not\in\{0,p-1\}$, then
  \begin{equation*}
    \begin{array}{c}
      \Ext_{U_\chi(W)}(V_\chi(\lambda'),V_\chi(\lambda))\\ \ \\
      \simeq \begin{cases}
        K,& \text{if }[\lambda-\lambda']\in\{2,3\},\\
        K,& \text{if }[\lambda-\lambda']=4 \text{ and }p\neq5,\\
        K,& \text{if }[\lambda-\lambda']=6 \text{ and }2\lambda=5\pm\sqrt{19}\text{ and }p>7,\\
        0,&\text{otherwise}.
      \end{cases}
    \end{array}
  \end{equation*}
\end{enumerate}
\end{Theorem}

\section{Height $-1$}\label{ch3}
Throughout this section we will assume that $r(\chi)=-1$, or equivalently, that $\chi=0$. Furthermore, we let as usual $\lambda$ and $\lambda'$ be elements in $\Lambda(0)\simeq\F_p$.
\subsection{$V_0(\lambda)$ and $V_0(\lambda')$}
This case is fully described in Section \ref{Verma}. The reader is referred to Theorem \ref{sum1}.
\subsection{$V_0(\lambda)$ and $K$}
Our approach will follow the one taken earlier in Section
\ref{Verma}. We have an isomorphism of vector spaces
\[
\Ext_{U_0(W)}(V_0(\lambda),K)\simeq \Ext_{U_0(W_0)}(K_{\lambda},K),
\]
where $K$ denotes the trivial $W$-module. Suppose that we have a short exact sequence of $U_0(W_0)$-modules
\begin{equation}\label{seq3}
0\rightarrow K\xrightarrow{f} M\xrightarrow{g} K_\lambda\rightarrow 0,
\end{equation}
and fix bases $v$ and $v'$ for $K$ and $K_\lambda$, respectively.  Let $\{w,w'\}$ be a basis for $M$ such that $f(v)=w$ and $g(w')=v'$. Since $U_0(Ke_0)$ is semisimple we may choose $w'$ such that $e_0w'=\lambda w'$. Note that the $\lambda$ weight space in $M$ is one dimensional so $w'$ is unique with these properties. Now, we have clearly $e_iw=0$ for all $i$. Furthermore, if $e_1w'$ and $e_2w'$ are both 0, then $e_iw'=0$ for all $i>0$. This can be seen by a simple induction argument since
\[
(i-1)e_{i+1}w'=(e_1e_i-e_ie_1)w'.
\]
The above observations give a necessary condition for \eqref{seq3} to be non-split, namely, that $e_1w'\neq0$ or $e_2w'\neq0$. We can interpret this condition in terms of $\lambda$. Indeed, since $e_iw'$ belongs to the weight space $M_{\lambda+i}$, we have $e_iw'=0$ if $\lambda+i\neq0$ and $i>0$. Thus $e_1w'=0$ and $e_2w'=0$ for $\lambda\not\in\{p-1,p-2\}$. This leads to the following lemma  
\begin{Lemma} If $\lambda\not\in\{p-1,p-2\}$ then
\begin{equation*} 
\Ext_{U_0(W_0)}(K_\lambda,K)=0.
\end{equation*}
\end{Lemma}
We consider the cases $\lambda\in\{p-1,p-2\}$ separately. Note that $\lambda=p-1$ is not interesting for our purpose because $V_0(p-1)$ is not simple; we shall, nevertheless, include it in our study for the sake of completeness. Fix $\lambda=p-k$ for $k\in\{1,2\}$. There exists $\a \in K$ such that for every $i>0$
\[
e_iw'=
\begin{cases}
\a w, &\text{if } i=k,\\
0, &\text{otherwise}.
\end{cases}
\]
As in Section \ref{Verma}, we construct a well-defined injective homomorphism
\[
\Ext_{U_0(W_0)}(K_\lambda,K)\rightarrow K,
\] 
which sends the class of \eqref{seq3} to $\a$. Conversely, for each $\a\in K$ consider a short exact sequence of vector spaces
\begin{equation*}
0\rightarrow K\xrightarrow{f} M_{\a}\xrightarrow{g} K_\lambda\rightarrow 0,
\end{equation*}
and choose a basis $\{w,w'\}$ for $M_{\a}$ such that $f(v)=w$ and $g(w')=v'$ where, as before, $v$ and $v'$ are bases for $K$ and $K_\lambda$, respectively. For each  $0\leq i\leq p-2$ define an endomorphism $E_i\in \End_K(M_{\a})$ of $M_{\a}$ such that $E_iw=0$ for all $i$ and 
\[
E_iw'=
\begin{cases}
\lambda w',&\text{if }i=0,\\
\a w,&\text{if }i=k,\\
0,&\text{otherwise}.
\end{cases}
\]
Furthermore, set $E_i^{[p]}=\delta_{i0}E_i$ for all $i=0,1,\dots,p-2$
and $E_i=0$ for $i\neq0,1,\dots,p-2$. We have
\[
(E_i^p-E_i^{[p]})w'=0\quad \text{for all} \quad 0\leq i\leq p-2,
\]
and
\[
[E_i,E_j]w'=(j-i)E_{i+j}w' \ \text{ for all}\ i \text{ and }j.
\]
To prove the above formula one can assume $i\in\{0,k\}$ and that
$i\neq j$ since otherwise one would obtain 0 on both sides of the equation.
We summarize the results in this section into the following Proposition
\begin{Proposition}We have
\[
\Ext_{U_0(W)}(V_0(\lambda),K)\simeq
\begin{cases}
K,&\text{if }\lambda\in\{p-1,p-2\},\\
0,&\text{otherwise}.
\end{cases}
\]
\end{Proposition}
Evidently, $K$ is a self-dual module in the sense that it is
isomorphic to its dual. Thanks to the isomorphism 
\[
\Ext_{U_0(W)}(K,V_0(\lambda))\simeq \Ext_{U_0(W)}(V_0(\lambda)^*,K^*),
\]
and the fact that $V_0(\lambda)^*$ and $V_0(p-1-\lambda)$ are isomorphic we obtain 
\begin{Corollary} We have
\[
\Ext_{U_0(W)}(K,V_0(\lambda))\simeq
\begin{cases}
K,&\text{if }\lambda=\{0,1\},\\
0,&\text{otherwise}.
\end{cases}
\]
\end{Corollary}

\subsection{$V_0(\lambda)$ and $S$}\label{3.1}We proceed to describe
$\Ext_{U_0(W)}(V_0(\lambda),S)$ where $S$ denotes the
$(p-1)$-dimensional simple $U_0(W)$-module. We have an isomorphism 
\[
\Ext_{U_0(W)}(V_0(\lambda),S)\simeq \Ext_{U_0(W_0)}(K_{\lambda},S).
\]
Thus, classifying the extensions of $V_0(\lambda)$ by $S$ reduces to classifying the extensions of $K_\lambda$ by $S$. Suppose that we have a short exact sequence of $U_0(W_0)$-modules
\begin{equation}\label{seq4}
0\rightarrow S\xrightarrow{f} M\xrightarrow{g} K_\lambda\rightarrow 0,
\end{equation}
and fix bases $v_0,v_1,\dots,v_{p-2}$ and $v'$ of $S$ and $K_\lambda$, respectively.  Let 
\begin{equation}\label{b1}
w_0,\dots,w_{p-2},w'
\end{equation}
 be a basis of $M$ such that $f(v_i)=w_i$ for all $i$ and $g(w')=v'$; we may choose $w'$ such that $e_0w'=\lambda w'$. We have 
\begin{equation*}
  e_jw_i=
  \begin{cases}
    (-1)^{j+1}\frac{(i+1)!}{(i-j)!}w_{i-j}, & \text{if }j\leq i,\\
    0, & \text{otherwise}.
  \end{cases}
\end{equation*}
By weight considerations there exists $\a_j\in K$ ($j>0$) such that
\[
e_jw'=\begin{cases}
\a_jw_{[-\lambda-j-1]},&\text{if }\lambda+j\neq0,\\
0,&\text{otherwise}.
\end{cases}
\]
Furthermore, we set 
\begin{equation}\label{addcond}
\a_j=0\quad\text{if}\quad\lambda+j=0 \ \text{ or }\  j\not\in\{1,2,\dots,p-2\},
\end{equation}
For $\lambda\neq0$ the $\lambda$ weight space in $M$ is two-dimensional; it is generated by $w'$ and $w_{[-\lambda-1]}$ so any different choice for $w'$ has the form $w'+bw_{[-\lambda-1]}$ for some
$b\in K$. Obviously, this leads to the same $e_jw'$ if
$e_jw_{[-\lambda-1]}=0$. In particular, we see that all the $\mathfrak{a}_j$ are determined by $M$ if $\lambda=p-1$. The same holds for $\lambda=0$ since $w'$ is unique in this case; we have $M_0=Kw'$.
\begin{Lemma}[A]If $\lambda\in\{0,p-1\}$ then all the $\a_j$ are determined by $M$.
\end{Lemma}
For $\lambda\not\in\{0,p-1\}$ we have $e_1w_{[-\lambda-1]}\neq0$. Thus, the discussion preceding Lemma A shows that we can choose $w'$ uniquely such that $e_1w'=0$.
\begin{Lemma}[B]
If $\lambda\not\in\{0,p-1\}$ then there is a unique choice for $w'$ such that $e_1w'=0$.
\end{Lemma} 
We will henceforth always assume that $w'$ is the unique choice from
Lemma B if $\lambda\not\in\{0,p-1\}$. As in Section \ref{Verma} we
obtain an injective homomorphism
$\Ext_{U_0(W_0)}(K_{\lambda},S)\rightarrow K^2$ which sends the class
of \eqref{seq4} to $(\a_1,\a_2)\in K^2$ and such that $\a_1=0$ for
$\lambda\not\in\{0,p-1\}$. We have
\begin{equation}\label{extend}
    \a_j=A_j\a_1+B_j\a_2 \quad\text{for all}\quad 3\leq j\leq p-2,
  \end{equation}
  where
  \begin{align*}
    A_j&=-\frac{1}{j-2}\biggl(\prod_{k=1}^{j}(\lambda+k)\biggr)\biggl(1+\sum_{k=4}^{j}\frac{(j-k)!}{(j-3)!}\prod_{l=1}^{k-3}(\lambda+j-l)\biggr),
  \end{align*}
  (the summation $\sum_{k=4}^j$ is understood to be 0 when $j=3$) and
  \[
  B_j=\frac{(-1)^j}{(j-2)!}\prod_{k=2}^{j-1}(-\lambda-1-k)(\lambda+k).
  \]
Conversely, let $(\a_1,\a_2)$ be a pair in $K^2$ such that, cf. \eqref{addcond}
\[
\a_1=0\quad\text{if}\quad\lambda=p-1,
\]
\[
\a_2=0\quad\text{if}\quad\lambda=p-2.
\]
We extend $(\a_1,\a_2)$ to a tuple $\a=(\a_1,\a_2,\dots,\a_{p-2}) \in K^{p-2}$ by using \eqref{extend} and consider a short exact sequence
\begin{equation}\label{seq6}
0\rightarrow S\xrightarrow{f} M_\a\xrightarrow{g} K_\lambda\rightarrow 0,
\end{equation}
where $M_\a$ is a vector space with a basis as in \eqref{b1}. For each $0\leq j\leq p-2$ we define an endomorphism $E_j\in \End_K(M_\a)$ such that 
\begin{equation*}
  E_jw_i=
  \begin{cases}
    (-1)^{j+1}\frac{(i+1)!}{(i-j)!}w_{i-j}, & \text{if }j\leq i,\\
    0, & \text{otherwise},
  \end{cases}
\end{equation*}
and $E_0w'=\lambda w'$ and for every $1\leq j\leq p-2$
\[
E_jw'=\begin{cases}
\a_jw_{[-\lambda-j-1]},&\text{if }\lambda+j\neq0,\\
0,&\text{otherwise}.
\end{cases}
\]
Furthermore, we set $E_j^{[p]}=\delta_{j0}E_j$ for all $j$ and $E_j=0$ for $j\not\in\{0,1,\dots,p-2\}$. We have
\[
(E_j^p-E_j^{[p]})w'=0 \quad \text{for all}\quad j.
\]
Except for the change in notation, the formula
$[E_i,E_j]w'=(j-i)E_{i+j}w' $ leads to the same equations as in
Section \ref{W--module}; there is no need to include them here.
\begin{Proposition} We have 
\[
\Ext_{U_0(W)}(V_0(\lambda),S)\simeq
\begin{cases}
K,&\text{if }\lambda\in\{0\}\cup\{p-i\mid 3\leq i\leq5\},\\
0,&\text{otherwise}.
\end{cases}
\]
\end{Proposition}
\begin{proof} The short exact sequence 
\begin{equation}\label{eq:5}
0\rightarrow S\rightarrow V_0(0)\xrightarrow{\pi} K\rightarrow 0
\end{equation}
induces the long exact sequence of vector spaces
\begin{align}\label{seq7}
0&\rightarrow\Hom_{U_0(W)}(V_0(\lambda),S)\rightarrow\Hom_{U_0(W)}(V_0(\lambda),V_0(0)) \\ \notag 
&\rightarrow\Hom_{U_0(W)}(V_0(\lambda),K)\rightarrow\Ext_{U_0(W)}(V_0(\lambda),S)\\ \notag
&\rightarrow\Ext_{U_0(W)}(V_0(\lambda),V_0(0))\rightarrow\Ext_{U_0(W)}(V_0(\lambda),K)\rightarrow\cdots.
\end{align}
Due to Schur's Lemma, we have $\Hom_{U_0(W)}(V_0(\lambda),K)=0$ for $\lambda\not\in\{0,p-1\}$.  We obtain, in these cases, the exact sequence 
\begin{align}\label{seq8}
0&\rightarrow\Ext_{U_0(W)}(V_0(\lambda),S)\rightarrow\Ext_{U_0(W)}(V_0(\lambda),V_0(0))&\\ \notag
&\rightarrow\Ext_{U_0(W)}(V_0(\lambda),K)\rightarrow\cdots.
\end{align}
Now $\Hom_{U_0(W)}(V_0(0),K)$ is 1-dimensional and generated by the surjection $\pi$ in \eqref{eq:5}. The map $\Hom_{U_0(W)}(V_0(0),V_0(0))\rightarrow\Hom_{U_0(W)}(V_0(0),K)$ appearing in \eqref{seq7} is surjective because it maps the identity $\id_{U_0(W)}$ on $U_0(W)$ to $\pi$. Thus, we obtain an exact sequence as in \eqref{seq8}
\begin{align*}
0&\rightarrow\Ext_{U_0(W)}(V_0(0),S)\rightarrow\Ext_{U_0(W)}(V_0(0),V_0(0))&\\
&\rightarrow\Ext_{U_0(W)}(V_0(0),K)\rightarrow\cdots.
\end{align*}
Consequently, we have for every $\lambda\neq p-1,p-2$
\[
\Ext_{U_0(W)}(V_0(\lambda),S)\simeq\Ext_{U_0(W)}(V_0(\lambda),V_0(0)),
\]
which, by Proposition \ref{sum1}, implies 
\[
\Ext_{U_0(W)}(V_0(\lambda),S)\simeq
\begin{cases}
K,&\text{if }\lambda\in\{0\}\cup\{p-i\mid 3\leq i\leq5\},\\
0,&\text{if }\lambda\not\in\{0\}\cup\{p-i\mid 1\leq i\leq5\}.
\end{cases}
\]
Next, suppose that $\lambda\in\{p-1,p-2\}$. We prove that every short
exact sequence of $U_0(W)$-modules as in \eqref{seq6} splits, or
equivalently, that $\a_1=0$ and $\a_2=0$. The case $\lambda=p-2$ comes
almost at once from the observations preceding the proposition; indeed
$\a_1=0$ because of Lemma B and $\a_2=0$ by definition, see
\eqref{addcond}. For $\lambda=p-1$ it follows from \eqref{addcond}
that $\a_1=0$. To prove $\a_2=0$ note first that \eqref{extend}
implies $\a_2=-2 \a_{p-2}$. However, since
\[
1\leq[-1-\lambda-(p-2)]\quad\text{and}\quad p-2\leq[-1-\lambda-1],
\]
we see by inserting $(i,j)=(1,p-2)$ into \eqref{A4} that $\a_{p-2}=0$. This is only possible if $\a_2=0$ so the proposition is proved.
\end{proof}
The module $S$ is a self-dual module. Thanks to the isomorphism
\[
\Ext_{U_0(W)}(S,V_0(\lambda))\simeq \Ext_{U_0(W)}(V_0(\lambda)^*,S^*)
\]
and the fact that  $V_0(\lambda)^*$ and $V_0(p-1-\lambda)$ are isomorphic we obtain
\begin{Corollary} We have
\begin{equation*}
\Ext_{U_0(W)}(S,V_0(\lambda))\simeq
\begin{cases}
K,&\text{if }\lambda\in\{2,3,4,p-1\},\\
0,&\text{otherwise}.
\end{cases}
\end{equation*}
\end{Corollary}
\subsection{Self-extensions of $K$}
There are no nontrivial self-extensions of the trivial module $K$.
\begin{Proposition} We have
\begin{equation*}
\Ext_{U_0(W)}(K,K)=0
\end{equation*}
\end{Proposition}
\begin{proof}
Suppose we have a short exact sequence of $U_0(W)$-modules
\begin{equation}\label{seq5}
0\rightarrow K\xrightarrow{f} M\xrightarrow{g} K\rightarrow 0,
\end{equation}
and let $\{w,w'\}$ be a basis of $M$ such that $Kw=\Img f$. We have clearly $Ww=0$. Furthermore, $W w'\subset Kw$ since $Ww'\subset\ker g$. Thus, every $x\in W$ acting on $M$ can be represented by a matrix 
\begin{equation*}
\begin{bmatrix}
0&\phi(x)\\
0&0
\end{bmatrix},
\end{equation*}
where $\phi:W\rightarrow K$ is a homomorphism of Lie algebras. That $\phi$ preserves the Lie algebra structure  follows from the fact that $M$ is a $W$-module and hence that $W\rightarrow \gl(M)\simeq\gl(2,K)$ is a homomorphism of Lie algebras. Now, if \eqref{seq5} is non-split, then $\phi(x)\neq0$ for some $x\in W$. The kernel of $\phi$ is then an ideal of codimension 1, in apparent contradiction with the fact that $W$ is simple.
\end{proof}
\subsection{Self-extensions of $S$}
There are no nontrivial self-extensions of the module $S$.
\begin{Proposition}We have
\[
\Ext_{U_0(W)}(S,S)=0.
\]
\end{Proposition}
\begin{proof}
The short exact sequence
\[
0\rightarrow S\rightarrow V_0(0)\rightarrow K\rightarrow 0
\]
induces the long exact sequence
\begin{align*}
0&\rightarrow\Hom_{U_0(W)}(S,S)\rightarrow\Hom_{U_0(W)}(S,V_0(0))\\
&\rightarrow\Hom_{U_0(W)}(S,K)\rightarrow\Ext_{U_0(W)}(S,S)\\
&\rightarrow\Ext_{U_0(W)}(S,V_0(0))\rightarrow\cdots,
\end{align*}
which by Corollary \ref{3.1} and Schur's lemma implies the claim.
\end{proof}
\subsection{$S$ and $K$}The proof of the next proposition makes
use of the fact that $\Hom_{U_0(W)}(V_0(0),S)=0$.
\begin{Proposition}We have
\[
\Ext_{U_0(W)}(K,S)\simeq K^2.
\]
\end{Proposition}
\begin{proof}
The short exact sequence
\[
0\rightarrow S\rightarrow V_0(0)\rightarrow K\rightarrow 0
\]
induces the long exact sequence
\begin{align*}
0&\rightarrow\Hom_{U_0(W)}(K,S)\rightarrow\Hom_{U_0(W)}(V_0(0),S)\\
&\rightarrow\Hom_{U_0(W)}(S,S)\rightarrow\Ext_{U_0(W)}(K,S)\\
&\rightarrow\Ext_{U_0(W)}(V_0(0),S)\rightarrow\Ext_{U_0(W)}(S,S)\rightarrow\cdots,
\end{align*}
which in turn induces the exact sequence
\begin{align*}
0&\rightarrow\Hom_{U_0(W)}(S,S)\rightarrow\Ext_{U_0(W)}(K,S)\rightarrow\Ext_{U_0(W)}(V_0(0),S)\rightarrow0.
\end{align*}
The claim follows from Schur's lemma and Proposition \ref{3.1}.
\end{proof}
Since $K$ and $S$ are self-dual modules, the isomorphism 
\[
\Ext_{U_0(W)}(K,S)\simeq \Ext_{U_0(W)}(S^*,K^*)
\]
implies the following corollary
\begin{Corollary} We have
\[
\Ext_{U_0(W)}(S,K)\simeq K^2.
\]
\end{Corollary}
\section{Height 0 and 1}
All the work has been done in Section \ref{Verma}; Theorem \ref{sum1}
gives a complete classification of the $\chi$-reduced Verma modules
having character $\chi$ of height 0 and 1. 
\section{Height $p-1$}
Recall that every simple $U_\chi(W)$--module with one exception is
projective. The remaining simple module $L$ has a projective cover
with two composition factors both isomorphic to $L$. We have a
short exact sequence
\[
0\rightarrow L\rightarrow P\rightarrow L\rightarrow 0,
\]
where $P$ is a projective module. This induces the long exact sequence
\begin{align*}
0&\rightarrow\Hom_{U_\chi(W)}(L,L)\rightarrow\Hom_{U_\chi (W)}(P,L)\\
&\rightarrow\Hom_{U_\chi (W)}(L,L)\rightarrow\Ext_{U_\chi (W)}(L,L)\\
&\rightarrow\Ext_{U_\chi (W)}(P,L)\rightarrow\cdots,
\end{align*}
which in turn implies $\Ext_{U_\chi (W)}(L,L)\simeq K$. We have therefore
\begin{Proposition} If $M$ and $N$ are two simple $U_\chi(W)$--modules then
\[
\Ext_{U_\chi(W)}(M,N)\simeq
\begin{cases}
K,&\text{if }N=M=L,\\
0,&\text{otherwise}.
\end{cases}
\]
\end{Proposition}
\bibliography{myBib}{}
\bibliographystyle{plain}
\end{document}